\documentclass[11pt]{article}

\newcommand{\documentdate}{6 IV 2021}

\usepackage{a4wide,latexsym,varioref,amsmath,amssymb,bbm,amsfonts,xcolor}

\newcommand{\numsection}[1]{\section{#1}\setcounter{equation}{0}}
\newcommand{\appnumsection}[1]{\section*{#1}\setcounter{equation}{0}
  \renewcommand{\theequation}{A.\arabic{equation}}
  \renewcommand{\thetheorem}{A.\arabic{theorem}}
  \renewcommand{\thetable}{A.\arabic{table}}
  \renewcommand{\thefigure}{A.\arabic{figure}}
  \renewcommand{\thesection}{A} }
\renewcommand{\theequation}{\arabic{section}.\arabic{equation}}
\renewcommand{\thetable}{\arabic{section}.\arabic{table}}
\renewcommand{\thefigure}{\arabic{section}.\arabic{figure}}

\newcommand{\calO}{{\cal O}}
\newcommand{\calS}{{\cal S}}

\newcommand{\beqn}[1]{\begin{equation}\label{#1}}
\newcommand{\eeqn}{\end{equation}}
\newcommand{\ii}[1]{\{1, \ldots, #1 \}}
\newcommand{\iiz}[1]{\{0, \ldots, #1 \}}
\newcommand{\iibe}[2]{\{ #1, \ldots, #2 \}}
\newcommand{\flow}{f_{\rm low}}
\newcommand{\barphi}{\overline{\phi}}

\newcommand{\DT}{\Delta T}
\newcommand{\barDT}{\overline{\Delta T}}

\newcommand{\sfrac}[2]{{\scriptstyle \frac{#1}{#2}}}
\newcommand{\half}{\sfrac{1}{2}}
\newcommand{\quarter}{\sfrac{1}{4}}
\newcommand{\tim}[1]{\;\; \mbox{#1} \;\;}
\newcommand{\req}[1]{(\ref{#1})}
\newcommand{\eqdef}{\stackrel{\rm def}{=}}
\newcommand{\ms}{\;\;\;\;}

\newcommand{\barf}{\overline{f}}

\newcommand{\barT}{\overline{T}}
\newcommand{\bigfrac}[2]{\frac{\displaystyle #1}{\displaystyle #2}}
\newcommand{\bigsum}{\displaystyle \sum}

\newtheorem{theorem}{Theorem}[section]
\newtheorem{lemma}[theorem]{Lemma}

\newcommand{\llem}[2]{\vspace{\baselineskip} 
\noindent\framebox[\textwidth]{\parbox{0.95\textwidth}{
\begin{lemma} \label{#1} \rm #2 \end{lemma} } } \vspace{\baselineskip} }
\newcommand{\lthm}[2]{\vspace{\baselineskip} 
\noindent\framebox[\textwidth]{\parbox{0.95\textwidth}{
\begin{theorem} \label{#1} \rm #2 \end{theorem} } } \vspace{\baselineskip} }
\newcommand{\bpr}{{\bf Proof.} \hspace{1.5mm}}
\newcommand{\epr}{\hfill $\Box$ \vspace*{1em}}
\newcommand{\proof}[1]{
\begin{list}{}{
\setlength{\topsep}{0.0pt}
\setlength{\partopsep}{0.0pt}
\setlength{\leftmargin}{0.025\textwidth}
\setlength{\rightmargin}{0.5\leftmargin}
\setlength{\labelwidth}{0.5\leftmargin}
\setlength{\labelsep}{0.25\leftmargin}}
\item \bpr #1 \epr \noindent
\end{list}}
\newcounter{algo}[section]
\renewcommand{\thealgo}{\thesection.\arabic{algo}}
\newcommand{\algo}[3]{\refstepcounter{algo}
\begin{center}\begin{figure}[htbp]
\framebox[\textwidth]{
\parbox{0.95\textwidth} {\vspace{\topsep}
{\bf Algorithm \thealgo : #2}\label{#1}\\
\vspace*{-\topsep} \mbox{ }\\
{#3} \vspace{\topsep} }}
\end{figure}\end{center}}
\newcommand{\al}[1]{{\footnotesize{\sf #1}}}
\newcommand{\accuracy}{{\tt accuracy}}
\newcommand{\accuracys}{{\tt accuracy}$_s$\,}
\newcommand{\accuracyj}{{\tt accuracy}$_j$\,}
\newcommand{\insufficient}{{\tt insufficient}}
\newcommand{\terminal}{{\tt terminal}}
\newcommand{\absolute}{{\tt absolute}}
\newcommand{\relative}{{\tt relative}}
\newcommand{\status}{{\tt status}}
\newcommand{\order}{{\tt order}}
\newcommand{\dlta}{{\tt delta}}
\newcommand{\orad}{{\tt radius}}
\newcommand{\innoisephi}{{\tt in-noise-phi}}
\newcommand{\innoises}{{\tt in-noise-s}}
\newcommand{\innoisef}{{\tt in-noise-f}}
\newcommand{\asol}{{\tt approximate-minimizer}}
\newcommand{\sphi}{\widehat{\phi}}
\newcommand{\epsmin}{\epsilon_{\min}}
\renewcommand{\Re}{\mathbb{R}}

\title{The Impact of Noise on Evaluation Complexity: \\ The Deterministic Trust-Region Case}

\author{
  S. Bellavia\thanks{Dipartimento di Ingegneria Industriale,
    Universit\`{a} degli Studi di Firenze, Italy. Member of the INdAM Research
    Group GNCS. Email: stefania.bellavia@unifi.it},
  G. Gurioli\thanks{Dipartimento di Matematica e Informatica ``Ulisse Dini'',
    Universit\`{a} degli Studi di Firenze, Italy.  {Member of the INdAM Research
    Group GNCS.} Email: gianmarco.gurioli@unifi.it},
  B. Morini\thanks{Dipartimento di Ingegneria Industriale,
    Universit\`{a} degli Studi di Firenze, Italy. Member of the INdAM Research
    Group GNCS. Email: benedetta.morini@unifi.it} \ and
  Ph. L. Toint\thanks{ Namur Center for Complex Systems (naXys),
    University of Namur, 61, rue de Bruxelles, B-5000 Namur, Belgium.
    Email: philippe.toint@unamur.be} 
}

\date{\documentdate}

\begin{document}


\maketitle

\begin{abstract}
Intrinsic noise in objective function and derivatives evaluations may cause
premature termination of optimization algorithms. Evaluation complexity bounds
taking this situation into account are presented in the framework of a
deterministic trust-region method. The results show that the presence of
intrinsic noise may dominate these bounds, in contrast with what is known for
methods in which the inexactness in function and derivatives' evaluations is
fully controllable. Moreover, the new analysis provides estimates of the
optimality level achievable, should noise cause early termination. It finally
sheds some light on the impact of inexact computer arithmetic on evaluation
complexity.
\end{abstract}


\numsection{Introduction}

This paper attempts to answer a simple question: how does noise in function
values and derivatives affect evaluation complexity of smooth optimization?
While analysis has been produced to indicate how high accuracy can be
reached by optimization algorithms even in the presence of inexact but
deterministic\footnote{Similar results are also known for the stochastic case
(see \cite{BandScheVice14,ChenMeniSche18,PaquSche20,BlanCartMeniSche19,BellGuriMoriToin20,BellGuri21}),
which is outside the scope of this paper.}
function and derivatives' values (see 
\cite{Cart91,ConnGoulToin00,XuRoosMaho20,BellGuriMori21,YaoXuRoosMaho20,GratSimoToin20,CartGoulToin20c}),
these approaches crucially rely on the assumption
that the inexactness is controllable, in that it can be made arbitrarily small
if required so by the algorithm. But what happens in practical applications
where significant noise is intrinsic and can't be assumed away? How is the
evaluation complexity of the optimization algorithm altered?

To limit the scope of this analysis, we focus here on trust-region methods for
unconstrained problems, a well known class of algorithms (see
\cite{ConnGoulToin00} for an in-depth coverage and \cite{Yuan15} for a more
recent survey), whose complexity was first investigated in
\cite{GratSartToin08}. We choose to base our present developements on the existing
analysis of \cite{CartGoulToin20c}, where the evaluation complexity of
trust-region methods with explicit dynamic accuracy is presented. It is shown
in this paper that, under standard Lipschitz continuity assumptions, a variant
of the classical trust-region algorithm using derivatives of degree one to $q$
and allowing the control of inexactness in objective function and derivatives'
values will find a $q$-th order $\epsilon$- approximate minimizer of the objective function in
$\calO(\epsilon^{-(q+1)})$ evaluations of $f$ and its derivatives.

Our purpose in this paper is to extend these results to the case where such
favourable assumptions of the noise can no longer be made, in that evaluation
of $f$ or its derivatives may simply fail if the requested accuracy is too high.
In that case, the desired $\epsilon$ optimality may not be reachable, and
our minimization algorithm may be forced to terminate before approximate
convergence can be declared. The question then arises to establish not only an
upper bound on the number of evaluations for this event to occur, but also
bounds, if possible, on the level of optimality achieved at termination.
However, since noisy problems often occur in a context where even moderate
accuracy is expensive to obtain, we wish our algorithms to preserve the
ability of the methods described in \cite{CartGoulToin20c,BellGuriMoriToin20}
to dynamically adjust accuracy requests in the limits imposed by noise.

\vspace*{2mm}
\noindent
{\bf Contributions.} 
We will present a trust-region method
  allowing dynamic accuracy control whenever possible given the level of
  noise, and show that termination of this algorithm will occur in at most
  $\calO\Big(\min[\vartheta_f^{-1},\vartheta_d^{-1}\epsilon^{-(q+1)},\epsilon^{-(q+1)}]\Big)$
    evaluations, where $\vartheta_f$ and $\vartheta_d$ are the absolute noise levels in
    $f$ and its derivatives, respectively, $\epsilon$ is the (ideally)
    sought optimality threshold and
    $q\geq1$ is the sought optimality order. In addition, we will derive upper
    bounds on the value of optimality measures at termination that depend on
    $\vartheta_f$. To the best of our knowledge, these results are the first
    of their kind. Finally, a brief discussion will illustrate our results in
    the case where intrinsic noise is caused by computer arithmetic and
    round-off errors.

\noindent
Because our development heavily hinges on \cite{CartGoulToin20c},
repeating some material from this source is necessary to keep our argument
understandable. We have however done our best to limit this repetition as much
as possible, pushing some of it in an appendix when possible.

Even if the analysis presented below does not depend in any way on the choice
of the optimality order $q$, the authors are well aware that, while requests
for optimality of orders $q\in\{1,2\}$ lead to practical, implementable
algorithms, this may no longer be the case for $q>2$, at least for now. For high orders, the
methods discussed in the paper therefore constitute an ``idealized'' setting
(in which complicated subproblems can be approximately solved without affecting the
evaluation complexity) and thus indicate the limits of achievable results.

\vspace*{2mm}
\noindent
    {\bf Outline.} A first section briefly recalls the context and the notion
    of high-order approximate minimizers.  Section~\ref{section:trqdan} then
    presents a ``noise-aware'' inexact trust-region algorithm and its
    evaluation complexity. Brief conclusions and perspectives are finally
    presented in Section~\ref{section:conclusion}. 

\vspace*{2mm}
\noindent
{\bf Basic notations.}
Unless otherwise specified, $\|\cdot\|$ denotes the standard
Euclidean norm for vectors and matrices.  For a general symmetric tensor $S$
of order $p$, we define
\[
\|S\| \eqdef \max_{\|v\|=1}  | S [v]^p |
= \max_{\|v_1\|= \cdots= \|v_p\|=1} | S[v_1, \ldots, v_p] |
\]
the induced Euclidean norm. We also denote by $\nabla_x^j f(x)$ the $j$-th
order derivative tensor of $f$ evaluated at $x$ and note that such a tensor is
always symmetric for any $j\geq 2$. $\nabla_x^0 f(x)$ is a synonym for $f(x)$.
$\lfloor \alpha \rfloor$ denotes the largest integer not exceeding
$\alpha$. For symmetric matrices, $\lambda_{\min}[M]$ is the
leftmost eigenvalue of $M$.

\numsection{High-Order Taylor Decrements and High-Order Optimality}
\label{section:high-order}

Throughout this paper, we consider the unconstrained problem given by
\beqn{problem}
\min_{x \in \Re^n}  f(x),
\eeqn
where we assume that
the \emph{values of the objective function $f$ and its derivatives are computed
  inexactly and are subject to noise}. Inexact quantities will be denoted by
an overbar, so that $\barf(s)$ is an inexact value of $f(x)$ and
$\overline{\nabla_x^j f}(x)$ an inexact value of $\nabla_x^j f(x)$.
We will also assume that

\vspace*{1mm}
\noindent
{\bf AS.1:} the objective function $f$ is $q$ times continuously
differentiable in $\Re^n$, for some $q\geq 1$;

\vspace*{1mm}
\noindent
{\bf AS.2:} the first $q$ derivative tensors of $f$ are globally Lipschitz
    continuous, that is, for each $j\in\ii{q}$ there exist a constant
    $L_{f,j} \geq 0$ such that, for all $x,y$ in $\Re^n$,
    \[
    \| \nabla_x^j f(x) - \nabla_x^j f(y) \| \leq L_{f,j} \| x-y \|;
    \]
    
\vspace*{1mm}
\noindent
{\bf AS.3:} the objective function $f$ is bounded below by $\flow$ on $\Re^n$.

\vspace*{1mm}
\noindent
In what follows, we consider algorithms that are able to exploit all available
derivatives of $f$.  As in many minimization methods, we would like to build a
model of the objective function $f$ using the truncated Taylor expansions
(now of degree $j$ for $j\in\ii{q}$) given by
\beqn{taylor}
T_{f,j}(x,s) \eqdef f(x) + \sum_{\ell=1}^j \nabla_x^\ell f(x)[s]^\ell,
\eeqn
where $\nabla_x^\ell f(x)$ is a $\ell$-th order symmetric tensor and
$\nabla_x^\ell f(x)[s]^\ell$ is this tensor applied to $\ell$ copies of the
vector $s$. More specifically, we will be interested, at a given iterate
$x_k$, in finding a step $s\in \Re^n$ which makes the \emph{Taylor decrements}
\[
\Delta T_{f,j}(x_k,s) \eqdef f(x_k)- T_{f,j}(x_k,s) =
T_{f,j}(x_k,0)-T_{f,j}(x_k,s)
\]
large (note that $\Delta T_{f,j}(x,s)$ is independent of $f(x)$). When this is
possible, we anticipate from the approximating properties of the Taylor
expansion that some significant decrease is also possible in $f$.  Conversely,
if $\Delta T_{f,j}(x,s)$ cannot be made large in a neighbourhood of $x$, we
must be close to an approximate minimizer.  More formally, we define, for some
\emph{optimality radius} $\delta \in (0,1]$, the measure
\beqn{phi-def}
\phi_{f,j}^{\delta}(x) = \max_{\|d\|\leq \delta} \Delta T_{f,j}(x,d),
\eeqn
that is the maximal decrease in $T_{f,j}(x,d)$ achievable in a ball of radius $\delta$
centered at $x$.  We then define $x$ to be a $q$-th order
$(\epsilon,\delta)$-approximate minimizer (for some accuracy requests
$\epsilon\in (0,1]^q$) if and only if 
\beqn{approx-min}
\phi_{f,j}^\delta(x) \leq \epsilon_j \frac{\delta^j}{j!}
\tim{ for } j \in \ii{q},
\eeqn
(a vector $d$ solving the optimization problem defining
$\phi_{f,j}^{\delta}(x)$ in \req{phi-def} is called an \emph{optimality
displacement}).  In other words, a $q$-th order
$(\epsilon,\delta)$-approximate minimizer is a point from which no significant
decrease of the Taylor expansions of degree one to $q$ can be obtained in a
ball of optimality radius $\delta$. This notion is coherent with standard
optimality measures for low orders\footnote{It is easy to verify that,
irrespective of $\delta$, \req{approx-min} holds for $j=1$ if and only if
$\|\nabla_x^1f(x)\|\leq\epsilon_1$ and that, if $\|\nabla_x^1f(x)\|=0$,
$\lambda_{\min}[\nabla_x^2 f(x)] \geq -\epsilon_2$ if and only if
$\phi_{f,2}^\delta(x) \leq \epsilon_2$.} and has the advantage of being
well-defined and continuous in $x$ for every order.

Unfortunately, the exact values of $f(x)$ and $\nabla_x^\ell f(x)$ may be
unavailable, and we then face several difficulties.
The first is that we can't consider the optimality measure \req{phi-def}
anymore, but could replace it by the inexact variant
\beqn{barphi-def}
\barphi_{f,j}^{\delta}(x) = \max_{\|d\|\leq \delta} \barDT_{f,j}(x,d).
\eeqn
where
\[
\barDT_{f,j}(x,d) \eqdef \barT_{f,j}(x,0)-\barT_{f,j}(x,d)
\tim{ with }
\barT_{f,j}(x_k,s)
\eqdef \barf(x_k) + \sum_{\ell=1}^j \overline{\nabla_x^\ell f}(x_k)[s]^\ell.
\]
However, computing the exact global maximum in this definition may also be
too expensive, and we follow  \cite[Theorem~6.3.5]{ConnGoulToin00} and
\cite{CartGoulToin20c} in choosing to use the approximate version
given by $\barDT_{f,j}(x,d)$, where
\beqn{bbarphi-def}
\varsigma\,\barphi_{f,j}^{\delta}(x) \leq \barDT_{f,j}(x,d),
\eeqn
for some displacement $d$ such that $\|d\|\leq \delta$ and some constant
$\varsigma \in (0,1]$.  Note that \req{bbarphi-def} does not assume the
knowledge of the global maximizer or $\barphi_{f,j}^{\delta}(x)$, but merely
that we can ensure \req{bbarphi-def} (see
\cite{deKlLaur19,deKlLaur20,SlotLaur20} for research in this direction). Note
also that, by definition, 
\beqn{weak-opt}
\barDT_{f,j}(x,d) \leq \varsigma \alpha
\tim{ implies }
\barphi_{f,j}^{\delta}(x) \leq \alpha.
\eeqn

The second difficulty occurs when computing a step $s_k$ which is supposed to
make the exact Taylor decrement $\DT_{f,j}(x_k,s_k)$ large, 
since we now have to resort to making the inexact decrement
\[
\barDT_{f,j}(x,s_k) \eqdef \barT_{f,j}(x_k,0)-\barT_{f,j}(x_k,s_k)
\]
large. It is therefore necessary to ensure, somehow, that the error on this
decrement does not dominate its value. The theory developed in this paper 
depends on making the \emph{relative} error on $\barDT_{f,j}(x_k,s_k)$ (for a
chosen step $s_k$) smaller than one, which is to require that 
\beqn{relerr-barDT}
|\barDT_{f,j}(x_k,s_k) - \DT_{f,j}(x_k,s_k)|
\leq \omega \barDT_{f,j}(x_k,s_k)
\eeqn
for some constant $\omega\in (0,1)$ to be specified later.
It is clearly not obvious at this point how to enforce this relative error
bound. For now, we simply assume that it can be done in a finite
number of evaluations of $\{\overline{\nabla_x^{\ell}f}(x)\}_{\ell=1}^j$ which
are inexact approximations of $\{\nabla_x^{\ell}f(x)\}_{\ell=1}^j$.  The
third difficulty arises when assessing the performance of the computed step:
is the predicted decrease in objective function predicted by the (inexact)
decrement significant in view of the (absolute) noise level in computing
$\barf(x_k)$ and $\barf(x_k+s)$?  If not, the obtained decrease is dominated
by noise in $f$ and thus unreliable.  To avoid this, our algorithms
will attempt to require that
\beqn{abserr-f}
|\barf(x_k)-f(x_k)| \leq \omega \barDT_{f,j}(x_k,s_k)
\tim{and}
|\barf(x_k+s_k)-f(x_k+s_k)| \leq \omega \barDT_{f,j}(x_k,s_k),
\eeqn
where $\omega$ is the parameter occuring in \req{relerr-barDT}.
The fourth, and for our present purpose, most significant difficulty is that
achieving \req{relerr-barDT} or \req{abserr-f} may require an accuracy 
of $f$ and its derivatives which is not feasible for noisy problems, and we
will have to prematurely terminate the algorithm. In what follows, we make the
assumption that this situation may occur (and thus does occur in the worst
case) if, for some $x_k$ of interest and $j\in\ii{q}$, 
\beqn{noise-termination}
|\barf(x_k)-f(x_k)| \leq \vartheta_f
\tim{ ~~or~~}
\|\overline{\nabla_x^\ell f}(x_k)-\nabla_x^\ell f(x_k)\|\leq \vartheta_d
\tim{for some $\ell \in \ii{j}$.}
\eeqn
for some non-negative \emph{absolute noise levels}  $\vartheta_f$ and $\vartheta_d$.  
The rest of our analysis will therefore focus on analyzing trust-region
algorithms which ensure that
\req{relerr-barDT} and \req{abserr-f} hold  as long as \req{noise-termination}
fail.

Like many trust-region methods, our proposed algorithms will consist  of an
initialization followed by a loop, performed until termination, in which one
successively  
\begin{enumerate}
\vspace*{-2mm}
\item \label{s-term} evaluates the function's derivatives and checks for termination,
\vspace*{-2mm}
\item \label{s-step} computes a step $s_k$ which approximately minimizes an (inexact) Taylor
  model $\barT_{f,j}(x_k,s)$  while remaining the inequality $\|s_k\| \leq \Delta_k$,
  where $\Delta_k$ is the current trust-region radius,%
\vspace*{-2mm}
\item \label{s-compf} evaluates the objective function at the new potential iterate and
  accepts or refuses the step,
\vspace*{-2mm}
\item updates the trust-region radius $\Delta_k$.
\vspace*{-2mm}
\end{enumerate}
The discussion above suggests that, at the very least, specialized versions of
the first three steps will be necessary.

\numsection{A Trust-Region Algorithm with Explicit Dynamic Accuracy and Noise}
\label{section:trqdan}

Because our analysis is based on  \req{relerr-barDT}
and \req{abserr-f}, we have to discuss how these conditions can be achieved.
For this purpose, we will use
the ``Explicit Dynamic Accuracy'' (EDA) framework already used in
\cite{ConnGoulToin00,BellGuriMoriToin19,GratSimoToin20}, in which absolute
accuracies on the function and derivatives values may be specified by the
algorithm by imposing the bounds
\beqn{eda-f}
|\barf(x)-f(x)| \leq \zeta_f
\eeqn
and
\beqn{eda-d}
\|\overline{\nabla_x^\ell f}(x)- \nabla_x^\ell f(x)\| \leq \zeta_d
\tim{for} \ell\in\ii{j}
\eeqn
before the actual computation of $\barf(x)$ and $\overline{\nabla_x^\ell f}(x)$
takes place\footnote{We could obviously use values of $\zeta_d$ and
  $\vartheta_d$ depending on the degree $\ell$, but we prefer the above
  formulation to simplify notations.}.
Such a framework is applicable for instance to multiprecision
computations \cite{High17,GratToin20} or to problems where the desired values
are computed by an iterative process whose accuracy can be monitored. In our
trust-region algorithm, the thresholds $\zeta_f$ and $\zeta_d$ will be
adaptively updated in the course of the iterations, but it is already clear
that requesting $\zeta_d < \vartheta_d$  will be impossible when
 \req{noise-termination} holds.

\subsection{Checking the accuracy of the model decrease}

However, before this happens, the algorithm will need to verify that the model
decrease relative accuracy bound \req{relerr-barDT} holds when the
``derivative-by-derivative'' absolute errors $\zeta_d$ are
known. As it turns out, this request has to be relaxed somewhat whenever the
right-hand side $\omega \barDT_{f,j}(x_k,s_k)$ is small, as can be expected
near a minimizers, and we have to replace the relative accuracy bound
\req{relerr-barDT} by an absolute error bound in that case. The management of
these crucial details is the object of the \al{CHECK} algorithm
\vpageref{algorithm:check}. To describe this algorithm in a general context,
we suppose that we have a $r$-th degree Taylor series $T_r(x,v)$ of a given
function about $x$ in the direction $v$, along with an 
approximation $\barT_r(x,v)$ and its decrement $\barDT_r(x,v)$.  
Additionally, we suppose that a bound $\delta \geq \|v\|$ is given, and that
\emph{required} relative and absolute accuracies $\omega$ and $\xi>0$ 
are on hand. Moreover, we assume that the \emph{current}
upper bound $\zeta_{d,i_\zeta}$ on absolute accuracies of the derivatives of
$\overline{T}_r(x,v)$ with respect to $v$ at $v=0$ are provided.
Because it will always be the case when we need it, 
we will assume for simplicity that $\barDT_r(x,v) \geq 0$.  Finally, the
relative accuracy constant $\omega \in (0,1)$ will be fixed in our trust-region
algorithm, and we assume that it is given when needed in \al{CHECK}. The
constants $\gamma_\zeta$, $\vartheta_f$ and $\vartheta_d$ of \req{noise-termination} are also
assumed to be known.

\algo{algorithm:check}{The \al{CHECK} algorithm}{
\vspace*{-3mm}
\[
\accuracy = \mbox{\al{CHECK}}\Big(\delta,\barDT_r(x,v),\zeta_{d,i_\zeta},\xi\Big).
\]
\vspace*{-2mm}
\begin{description}
\item[\hspace*{6.5mm}] 
  If
  \vspace*{-3mm}
  \beqn{verif-term-2}
  \barDT_r(x,v) > 0
  \tim{ and }
  \zeta_{d,i_\zeta}\sum_{\ell=1}^r \frac{\delta^\ell}{\ell!} \leq \omega \,\barDT_r(x,v),
  \vspace*{-3mm}
  \eeqn
  set \accuracy\ to \relative.
\item[\hspace*{6.5mm}] 
  Otherwise, if
  \vspace*{-3mm}
  \beqn{verif-term-3}
  \zeta_{d,i_\zeta} \sum_{\ell=1}^r \frac{\delta^\ell}{\ell!} \leq \omega \,\xi\, \frac{\delta^r}{r!},
  \vspace*{-3mm}
  \eeqn
  set \accuracy\ to \absolute.
\item[\hspace*{6.5mm}] 
  Otherwise, if
  \beqn{improve-possible}
  \gamma_\zeta \zeta_{d,i_\zeta} > \vartheta_d,
  \eeqn
  set 
  \beqn{improve-accuracy}
  \zeta_{d,i_\zeta+1}= \gamma_\zeta \zeta_{d,i_\zeta}
  \eeqn
  and \accuracy\ to \insufficient.
\item[\hspace*{6.5mm}] 
 Otherwise, set \accuracy\ to \terminal.
\end{description}
}

\noindent
Note that the integer $i_\zeta$ counts the number of times the accuracy
threshold has been reduced by a factor $\gamma_\zeta$.
The outcome of the \al{CHECK} algorithm can then characterized as follows.

\llem{lemma:check}{
Let $\omega \in (0,1)$ and $\delta$, $\xi$ and $\zeta_{d,i_\zeta}$ be positive. 
Suppose that $\barDT_r(x,v) \geq 0$ and \req{eda-d} hold.  Then the call
$\accuracy = \mbox{\al{CHECK}}\Big(\delta,\barDT_r(x,v),\zeta_{d,i_\zeta},\xi\Big)
$
ensures that

\noindent
{\em (i)} \accuracy\ is either \absolute\ or \relative\ whenever
\vspace*{-3mm}
\[
\zeta_{d,i_\zeta}\sum_{\ell=1}^r \frac{\delta^i}{i!} \leq \omega \xi \frac{\delta^r}{r!};
\]
 
\noindent
{\em (ii)} if \accuracy\ is \absolute, 
\[
\max\Big[\barDT_r(x,v),
    \left|\barDT_r(x,w)- \Delta T_r(x,w)\right|\Big]
  \leq \xi \frac{\delta^r}{r!}
\]
for all $w$ with $\|w\|\leq \delta$;

\noindent
{\em (iii)} if \accuracy\ is \relative, 
  $\barDT_r(x,v)>0$ and
\[
  \left|\barDT_r(x,w)- \Delta T_r(x,w)\right|
  \leq \omega \barDT_r(x,v),
  \tim{for all $w$ with $\|w\|\leq \delta$}.
\]

\noindent
Moreover, the outcome $\accuracy = \insufficient$ indicates that new values of
the required approximate derivatives should be computed with the updated accuracy
thresholds, while $\accuracy = \terminal$ indicates that the minimization
algorithm has reached the noise level and should be terminated. 
}

\proof{We note that the \al{CHECK} algorithm is identical to the \al{VERIFY}
  algorithm of \cite{CartGoulToin20c} (itself inspired by
  \cite{BellGuriMoriToin19})  whenever \accuracy\ is either
  \absolute\ or \relative. Lemma~2.1 in that reference therefore ensures the
  conclusions (i) to (iii). If $\accuracy = \insufficient$, then
  \req{improve-possible} ensures that the accuracy threshold update
  \req{improve-accuracy} has been performed safely
  (\req{noise-termination} remains violated), while $\accuracy = \terminal$
  indicates that this was not the case, suggesting termination.
}

\noindent
Note that case (ii) is the case where relative accuracy would be excessively
requiring and absolute accuracy is preferred.

\subsection{Testing for termination}

We now start constructing  our new algorithm (which we call the
\al{TR$q$EDAN} algorithm because it uses the \al{EDA} framework and handles
\al{N}oise) step by step, following the outline given at the end of
Section~\ref{section:high-order}. Consider Step~1 first.
Since we have to rely on $\overline{\nabla_x^\ell f}(x_k)$ rather than
$\nabla_x^\ell f(x_k)$, it is clear that our optimality measure \req{phi-def} and test
\req{approx-min} should be modified to use the inexact values.
Ideally, we could mimic \cite{CartGoulToin20c} and terminate because of
\req{weak-opt} as soon as 
\beqn{approx-min-inexact}
\barDT_{f,j}(x_k,d_{k,j}) \leq \left(\frac{\varsigma\epsilon_j}{1+\omega}\right) \frac{\delta_k^j}{j!}
\tim{ for } j \in \ii{q},
\eeqn
and where $\omega\in(0,1)$ is the (still unspecified) relative accuracy
parameter of \req{relerr-barDT},
\[
\varsigma \barphi_{f,j}^{\delta_k}(x_k) \leq \barDT_{f,j}(x_k,d_{k,j})
\]
and $\delta_k$ is the optimality radius
at iteration $k$ (which we leave again unspecified at this stage).
However, we now have to take into account the
fact that noise in the values of the derivatives may prevent a meaningful
computation of $\barDT_{f,j}(x_k,d_{k,j})$. We therefore have to modify the
technique proposed in \cite[Algorithm~2.2]{CartGoulToin20c}.  Assuming that
the optimality radius $\delta_k$ is given, we thus consider
Algorithm~\ref{algorithm:phij} for computing the $j$-th approximate optimality measure
which is needed in \req{approx-min-inexact} to test for termination in the first step of the
\al{TR$q$EDAN} algorithm.

\algo{algorithm:phij}{Computing \boldmath$\barDT_{f,j}(x_k,d_{k,j})$}{ 
The iterate $x_k$, the index $j\in\ii{q}$ and the radius $\delta_k\in (0,1]$
are given, as well as constants $\gamma_\zeta \in (0,1)$ and
$\varsigma \in (0,1]$. The counter $i_\zeta$, the relative accuracy
$\omega\in (0,1)$ and the absolute accuracy bound $\zeta_{d,i_\zeta}$ are also given. 
\\
\vspace*{-3mm}
\begin{description}
\item[Step 1.1: ] If they are not yet available, compute
  $\{\overline{\nabla_x^if}(x_k)\}_{i=1}^j$ satisfying \req{eda-d} for
  $\zeta_d=\zeta_{d,i_\zeta}$.
\item[Step 1.2: ]
     Find $d_{k,j}$ with $\|d_{k,j}\|\leq \delta_k$ such that
     $
     \varsigma \barphi_{f,j}^{\delta_k}(x_k) \leq \barDT_{f,j}(x_k,d_{k,j})
     \vspace*{-1mm}
     $ 
     and compute
     \beqn{check-for-phi}
     \mbox{\accuracyj} =
     \mbox{\al{CHECK}}\Big(\delta_k,\barDT_{f,j}(x_k,d_{k,j}),
             \zeta_{d,i_\zeta}, \half \varsigma \epsilon_j\Big).
     \eeqn
\item[Step 1.3: ] If \accuracyj is \absolute\ or \relative, return
      $\barDT_{f,j}(x_k,d_{k,j})$.
\item[Step 1.4: ] If \accuracyj is \insufficient, return to Step~1.1 (with the
     tightened accuracy threshold $\zeta_{d,i_\zeta+1}$).
     Else (i.e.\ if \accuracyj is \terminal), terminate the \al{TR$q$EDAN}
     algorithm with $\tilde{x}=x_k$, \status\ = \innoisephi, \order = $j$
     and  \dlta = \orad\ = $\delta_k$.
\end{description}
}
  
\noindent
Note that, when termination occurs, this algorithm (and other algorithms we
will meet later) sets the four flags \status, \order, \dlta\ and \orad, which
will allow the user to determine the reason of termination once it occured
and, as we will see in Theorem~\ref{theorem:trqedan-complexity} below,
derive some useful properties of the point $\tilde{x}$ returned.

Because Algorithm~\ref{algorithm:phij} and
\cite[Algorithm~2.2]{CartGoulToin20c} only differ in Step~1.4, we may then
follow the reasoning of \cite[Lemma~2.2]{CartGoulToin20c} and obtain the
following result.

\llem{lemma:phi-res}{If Algorithm~\ref{algorithm:phij} terminates within Step~1.3 when 
\accuracyj is \absolute, then
\beqn{term-absolute}
\phi_{f,j}^{\delta_k}(x_k) \leq \epsilon_j \frac{\delta_k^j}{j!}.
\eeqn
Otherwise, if it terminates with \accuracyj being \relative, then
\beqn{term-relative}
(1-\omega) \barDT_{f,j}(x_k,d_{k,j})
\leq \phi_{f,j}^{\delta_k}(x_k)
\leq \left(\frac{1+\omega}{\varsigma}\right) \barDT_{f,j}(x_k,d_{k,j})
\eeqn
Moreover, termination with one of these two outcomes must occur if
\beqn{term-bound-phi}
\zeta_{d,i_\zeta} \leq \frac{\omega}{4}\,\varsigma\,\epsilon_j\,\frac{\delta_k^{j-1}}{j!}.
\eeqn
}

\noindent
Of course, termination may occur before \req{term-bound-phi} occurs (for
instance because of \req{noise-termination} 
in the call to \al{CHECK} in Step~1.2), but the bound
\req{term-bound-phi} shows that, if this doesn't happen, the accuracy
threshold $\zeta_{d,i_\zeta}$ can not be reduced infinitely often by the
factor $\gamma_\zeta$ and thus the loop between Steps~1.4 and 1.1 is
finite. Note that the rightmost inequality in \req{term-relative} and
\req{approx-min-inexact} together also imply \req{term-absolute} for order
$j$, justifying our choice of the scaling by $(1+\omega)$ in the former. 

Refering now to our outline on the trust-region method at the end of
Section~\ref{section:high-order}, we may now use
Algorithm~\ref{algorithm:phij} to implement a complete Step~1.  The idea is
first to identify a suitable optimality radius, which we choose to be
\beqn{deltak-def}
\delta_k = \min[\Delta_k, \theta]
\eeqn
(for some constant $\theta \leq 1$),
estimate the needed (inexact) derivatives and $\barphi_{f,j}^{\delta_k}(x_k)$  for
$j\in\ii{q}$ and decide on termination. The result is the \al{STEP1} algorithm
\vpageref{algorithm:STEP1}.

\algo{algorithm:STEP1}{\al{STEP1} for the \al{TR$q$EDAN} algorithm}{
  Set $\delta_k$ according to \req{deltak-def}. \\
  For $j = 1, \ldots, q $,\\
  \hspace*{2mm}\parbox[t]{14.5cm}{\vspace*{-2mm}\begin{description}
  \item{1.}  Evaluate $\overline{\nabla_x^j f}(x_k)$ and compute
  $\barDT_{f,j}(x_k,d_{k,j})$ using Algorithm~\ref{algorithm:phij}.
  \item{2.} If termination of the \al{TR$q$EDAN} algorithm has not
  happened in Step~1.4 of Algorithm~\ref{algorithm:phij} and
  \beqn{cauchy}
  \barDT_{f,j}(x_k,d_{k,j})
  > \left(\frac{\varsigma\epsilon_j}{1+\omega}\right) \frac{\delta_k^j}{j!},
  \eeqn
  exit \al{STEP1} with the current value of $j$ and the
  optimality displacement $d_{k,j}$ associated  with
  $\barphi_{f,j}^{\delta_k}(x_k)$. Otherwise consider the next $j$.
  \end{description}}\\*[1.5ex]
  Terminate the \al{TR$q$EDAN} algorithm with $\tilde{x} = x_k$,
  \status\ = \asol, \order\ = $q$ and \dlta\ = \orad\ =  $\delta_k$.
}

\noindent
Before progressing any further, we state an easy but useful technical
inequality.

\llem{lemma:chi-ineq}{
Consider $\delta\geq 0$.  Then, for all $j\geq 1$,
\beqn{chi-ineqs} 
\min[\delta,1] \leq \sum_{\ell=1}^j\frac{\delta^\ell}{\ell!} < 2
\max[\delta,\delta^j].
\eeqn
}

\proof{
  The bounds \req{chi-ineqs} easily follow from
  $
  1 \leq \bigsum_{\ell=1}^j\frac{1}{\ell!} < e-1 < 2.
  $
}

\noindent
We now consider what can be said if the \al{TR$q$EDAN} algorithm
terminates within \al{STEP1}.

\llem{lemma:step1-termination}{\mbox{}
\begin{itemize}
\item[(i)] Suppose that termination of the \al{TR$q$EDAN} algorithm occurs
  within \al{STEP1} with \status\ = \innoisephi, \order\ = $j$ and \dlta\ = $\delta_k$.
  Then
  \beqn{step1-termination}
  \phi_{f,i}^{\delta_k}(\tilde{x}) \leq \epsilon_i\frac{\delta_k^i}{i!}
  \tim{for} i\in\ii{j-1}
  \tim{ ~and~ }
  \phi_{f,j}^{\delta_k}(\tilde{x}) < \frac{4\vartheta_d}{\gamma_\zeta \omega}\delta_k.
  \eeqn 
\item[(ii)] Suppose that termination of the \al{TR$q$EDAN} algorithm occurs
  with \status\ = \asol\ and \dlta\ = $\delta_k$.
  Then \req{approx-min} holds and $\tilde{x}$ is a
  $q$-th order $(\epsilon,\delta_k)$-approximate minimizer.
\end{itemize}
}

\proof{
  We prove case (ii) first, which can only occur if Algorithm~\ref{algorithm:phij}
  terminates within Step~1.3 and \req{cauchy} fails for every $j\in\ii{q}$.
  We then have from Lemma~\ref{lemma:phi-res} that, for every $j\in\ii{q}$, 
  \[
  \phi_{f,j}^{\delta_k}(x_k)
  = \phi_{f,j}^{\delta_k}(\tilde{x})
  \leq \max\left[ \epsilon_j \frac{\delta_k^j}{j!},\left(\frac{1+\omega}{\varsigma}\right)\,\barDT_{f,j}(x_k,d_{k,j})\right]
  \leq \epsilon_j \frac{\delta_k^j}{j!},
  \]
  the last inequality resulting from the failure of \req{cauchy}.  Thus \req{approx-min} holds.

  Consider now case (i), that is when the call \al{CHECK} in Step~1.2 of
  Algorithm~\ref{algorithm:phij} returns \accuracyj $=$ \terminal\ for some $j\in\ii{q}$.
  Thus Algorithm~\ref{algorithm:phij} has terminated within Step~1.3 and
  \req{cauchy} has failed for every order of index smaller than
  ${j-1}$. Applying the same reasoning as for case (ii), we obtain that the
  first part of \req{step1-termination} holds. Now suppose that, instead of
  the call \req{check-for-phi} resulting in \accuracyj being \terminal, we had
  made the hypothetical call 
  \beqn{hypothetical-check}
  \mbox{\accuracyj} =
     \mbox{\al{CHECK}}\Big(\delta_k,\barDT_{f,j}(x_k,d_{k,j}),\zeta_{i,i_\zeta},
     \frac{\zeta_{d,i_\zeta}j!}{\omega\delta_k^j}\sum_{\ell=1}^j\frac{\delta_k^\ell}{\ell!}\Big).
  \eeqn
  Observe first that, since the call \req{check-for-phi} returned \terminal, \req{verif-term-2}
  failed on that call, and thus, since this is independent of the last
  argument of the call, it also fails for the call
  \req{hypothetical-check}. However, one easily checks that \req{verif-term-3}
  holds as an equality for this hypothetical call, and thus
  \req{hypothetical-check} would return \accuracyj\ as \absolute.  We may then
  use case (ii) in Lemma~\ref{lemma:check} and deduce from the triangular
  inequality that, for some $\tilde{d}$ with $\|\tilde{d}\|\leq \delta_k$,
  \[
  \phi_{f,j}^{\delta_k}(\tilde{x})
  = \Delta T_j(\tilde{x},\tilde{d})
  \leq \barDT_j(\tilde{x},\tilde{d})
       + \left|\barDT_j(\tilde{x},\tilde{d})- \Delta T_j(\tilde{x},\tilde{d})\right|
  \leq 2\,\frac{\zeta_{d,i_\zeta}j!}{\omega\delta_k^j}
       \left(\sum_{\ell=1}^j\frac{\delta_k^\ell}{\ell!}\right)\frac{\delta_k^j}{j!}.
  \]
  Moreover, since the call \req{check-for-phi} returned \terminal, we have
  that $\gamma_\zeta \zeta_{d,i_\zeta} < \vartheta_d$, and we deduce that
  \beqn{phi-terminal-j}
  \phi_{f,j}^{\delta_k}(\tilde{x})
  < 2 \, \frac{\vartheta_d}{\gamma_\zeta\omega}
       \left(\sum_{\ell=1}^j\frac{\delta_k^\ell}{\ell!}\right).
  \eeqn
  The second part of \req{step1-termination} then results from this inequality
  and \req{chi-ineqs} for $\delta=\delta_k\leq \theta\leq 1$.
} 

\noindent
We also have the following useful result.

\llem{lemma:no-term-step1}{
Suppose that, at iteration $k$, termination of the \al{TR$q$EDAN} algorithm
does not happen during execution of \al{STEP1}. Then
\beqn{barDTd-lower}
\barDT_{f,j}(x_k,d_{k,j})
\geq \frac{\zeta_{d,i_\zeta}}{\omega} \sum_{\ell=1}^j \frac{\delta_k^\ell}{\ell!},
\eeqn
where the threshold $\zeta_{d,i_\zeta}$ refers to its value at the end of
\al{STEP1}.
Moreover,
\beqn{phi-end-step1}
\phi_{f,i}^{\delta_k}(x_k) \leq \epsilon_i\frac{\delta_k^i}{i}
\tim{for} i\in \ii{j-1}
\tim{ and }
\phi_{f,j}^{\delta_k}(x_k) \leq \left(\frac{1+\omega}{\varsigma}\right)\barphi_{f,j}^{\delta_k}(x_k).
\eeqn
}

\proof{
Suppose that the last value of \accuracyj\ computed during the execution of
\al{STEP1} is \absolute.  Then Lemma~\ref{lemma:phi-res} implies that
\req{term-absolute} holds.  But, since $\omega \in (0,1)$, this  and
Lemma~\ref{lemma:check} (ii) contradict
\req{cauchy}. As a consequence, the last value of \accuracyj\ must be
\relative, in which case \req{verif-term-2} ensures \req{barDTd-lower}.
The first part of \req{phi-end-step1} again follows from the reasoning of
Lemma~\ref{lemma:step1-termination}(ii) for $i\in\ii{j-1}$. Finally,
the fact that \accuracyj\ is \relative\ implies that \req{term-relative} holds in
Lemma~\ref{lemma:phi-res}, which gives the second part of \req{phi-end-step1}.
}

\subsection{Computing a step}

Given Step~1, constructing Step~2 of our \al{TR$q$EDAN} algorithm is
relatively straightforward and we immediately provide the details in the
\al{STEP2} algorithm \vpageref{algorithm:STEP2}.

\algo{algorithm:STEP2}{\al{STEP2} for the \al{TR$q$EDAN} algorithm}{
The iterate $x_k$, the relative accuracy $\omega$, the requested accuracy
$\epsilon_j\in (0,1]^q$, the constant $\gamma_\zeta \in (0,1)$  the counter
$i_\zeta$ and the absolute accuracy threshold $\zeta_{d,i_\zeta}$ are given.
The index $j\in\ii{q}$, the optimality displacement $d_{k,j}$ resulting from
Step~1 and the constant $\theta \in (0,1]$, are also given such
that, by \req{cauchy}, 
\beqn{Cdecr}
\barDT_{f,j}(x_k,d_{k,j}) > \left(\frac{\varsigma\epsilon_j}{1+\omega}\right) \frac{\delta_k^j}{j!}.
\eeqn
\begin{description}
\item[Step 2.1: ] If they are not yet available, compute
    $\{\overline{\nabla_x^\ell f}(x_k)\}_{i=1}^j$ satisfying \req{eda-d} for
  $\zeta_d=\zeta_{d,i_\zeta}$ and recompute $\barDT_{f,j}(x_k,d_{k,j})$ for
  this accuracy threshold.
\item[Step 2.2: Step computation.]
  If $\Delta_k \leq \theta$, set $s_k = d_{k,j}$ and exit the \al{STEP2}
  algorithm with $\barDT_{f,j}(x_k,s_k) = \barDT_{f,j}(x_k,d_{k,j})$.
  Otherwise, find $s_k$ such that $\|s_k\|\leq\Delta_k$ and
  \beqn{step-subproblem}
  \barDT_{f,j}(x_k,s_k) \geq \barDT_{f,j}(x_k,d_{k,j}),
  \eeqn
  and compute
    \vspace*{-2mm}
    \beqn{check-sk}
    \mbox{\accuracys} = 
    \mbox{\al{CHECK}}\Big( \|s_k\|,\barDT_{f,j}(x_k,s_k),\zeta_{d,i_\zeta},
    \bigfrac{\varsigma\epsilon_j}{4(1+\omega)}\,\Big(\bigfrac{\theta}{\max\big[\theta,\|s_k\|\big]}\Big)^j\Big).
    \eeqn
  \item [Step 2.3:] If \accuracys\ is \relative, exit the \al{STEP2} algorithm
    with the step $s_k$ and the associated $\barDT_{f,j}(x_k,s_k)$.
\item[Step 2.4: ] If \accuracys is \insufficient, return to Step~2.1 (with the
  tightened accuracy thresholds).
  Otherwise, if \accuracys is \terminal, terminate the \al{TR$q$EDAN}
  algorithm  with $\tilde{x}=x_k$, \status\ = \innoises, \order\ = $j$,
  \dlta\ = $\delta_k$ and \orad\ = $\|s_k\|$.
\end{description}
}

  \noindent
Note that setting $s_k=d_{k,j}$ when $\Delta_k<\theta$ 
makes sense since $d_{k,j}$, computed in Step~1.2, is already a (\al{CHECK}ed)
approximate global maximizer of $\barDT_{f,j}(x_k,s)$ in the ball of radius
$\delta_k=\Delta_k$.
Two features of this algorithm remain nevertheless somewhat mysterious at this
stage.  The first is the complicated function of
$\|s_k\|$ and $\epsilon_j$ occuring in the last argument of the call to the
\al{CHECK} algorithm. As it turns out, it is possible to show
that the conjunction of \req{cauchy} and this particular
call to \al{CHECK}\footnote{\al{VERIFY} in \cite{CartGoulToin20c}.} ensures
that \accuracys\ cannot be \absolute.  This then also clarifies the second
mysterious feature, which is why this value of \accuracys\ is not considered in the rest
of the algorithm.  This is part of the following lemma, which was proved as
Lemma~3.2 in \cite{CartGoulToin20c} and which we can reuse since
the step computation
in that reference\footnote{In \cite{CartGoulToin20c}, the step computation is the combination
of Step~2 in Algorithm~3.1 and Algorithm~3.2 for the case where $\Delta_k\geq
\theta$. Note that, in this case, $\delta_k=\theta$ and thus $\delta_k$ may
be replaced by $\theta$ in the right-hand side of \req{Cdecr}, as
stated in Algorithm~3.2 of \cite{CartGoulToin20c}.}
and the \al{STEP2} algorithm only differ in the possibility that the \al{TR$q$EDAN}
algorithm can terminate in the call to \al{CHECK} in Step~2.2.

\llem{lemma:step2-res}{
  Suppose that the \al{TR$q$EDAN} algorithm does not terminate within the call
  to \al{CHECK} in Step~2.2 of the \al{STEP2} algorithm. Then the \al{STEP2}
  algorithm terminates with \accuracys being \relative\  and
  \req{relerr-barDT} holds. Moreover, this outcome must occur if 
  \beqn{step2-bound}
  \zeta_{d,i_\zeta}
  \leq \frac{\varsigma \omega\delta_k^j}{8j!(1+\omega)} \,\frac{\epsilon_j}{\max[1,\Delta_{\max}^j]}.
  \eeqn
}

\noindent
As for Lemma~\ref{lemma:phi-res}, the bound \req{step2-bound} ensures that the
loop between Steps~2.4 and 2.1 is finite. 

We conclude this paragraph by examining the optimality guarantees which may be
obtained, should the \al{TR$q$EDAN} algorithm terminate in \al{STEP2}.

\llem{lemma:step2-termination}
{
  Suppose that, at iteration $k$, the \al{TR$q$EDAN} algorithm terminates
  within  \al{STEP2} with \status\ = \innoises, \order = $j$  and \orad\ = $\|s_k\|$.
  Then
\beqn{step2-termination}
\phi_{f,j}^{\|s_k\|}(\tilde{x}) \leq \frac{4\vartheta_d}{\gamma_\zeta\omega}
\max\big[\|s_k\|,\|s_k\|^j\big].
\eeqn
}

\proof{
The fact that \status\ = \innoises\ implies that termination occurs in
Step~2.4, and it must be because the call \req{check-sk}
returns \accuracys\ equal to \terminal. As in the proof of
Lemma~\ref{lemma:step1-termination}, we consider replacing this call by the
hypothetical
\beqn{hypothetical-check-sk}
  \mbox{\accuracys} =
     \mbox{\al{CHECK}}\Big(\|s_k\|,\barDT_{f,j}(x_k,s_k),
 \zeta_{i,i_\zeta}, \frac{\zeta_{d,i_\zeta}j!}{\omega\|s_k\|^j}\sum_{\ell=1}^j\frac{\|s_k\|^\ell}{\ell!}\Big)
\eeqn
and verify that this call must return \accuracys\ equal to \absolute.  We also
deduce from case (ii) in Lemma~\ref{lemma:check}, the triangular inequality
and the bound $\gamma_\zeta\zeta_{d,i_\zeta}<\vartheta_d$ that, for some $\tilde{d}$
with $\|\tilde{d}\|\leq\|s_k\|$,
\[
\phi_{f,j}^{\|s_k\|}(\tilde{x})
= \Delta T_j(\tilde{x},\tilde{d})
\leq \barDT_j(\tilde{x},\tilde{d})
      + \left|\barDT_j(\tilde{x},\tilde{d})- \Delta T_j(\tilde{x},\tilde{d})\right|
\leq 2 \, \frac{\vartheta_d}{\gamma_\zeta\omega}\left(\sum_{\ell=1}^j\frac{\|s_k\|^\ell}{\ell!}\right),
\]
and \req{step2-termination} follows from \req{chi-ineqs}.
} 

\subsection{The complete \al{TR$q$EDAN} algorithm}

Having constructed the first two steps of the \al{TR$q$EDAN} algorithm, we are
now in position to specify the algorithm in its entirety
(see \vpageref{algorithm:trqedan}), making the necessary changes to handle
\req{noise-termination} in Step~3 along the way.  

\algo{algorithm:trqedan}{The \al{TR$q$EDAN} algorithm}
{
\vspace*{-2mm}
\begin{description}
\item[Step~0: Initialisation.]
  A criticality order $q$, a starting point $x_0$ and an initial trust-region radius
  $\Delta_0$ are given, as well as accuracy levels $\epsilon \in (0,1)^q$ and
  an  initial bound on absolute derivative accuracies $\kappa_\zeta$.
  The constants $\omega$, $\varsigma$, $\theta$, $\eta_1$, $\eta_2$,
  $\gamma_1$, $\gamma_2$,  $\gamma_3$ and $\Delta_{\max}$ are also given and satisfy
  \[
  \theta \in [\min_{j\in\ii{q}} \epsilon_j, 1],\;\;\;
  \Delta_0\leq \Delta_{\max},\;\;\;
  0 < \eta_1 \leq \eta_2 < 1, \;\;\;
  0< \gamma_1 < 1 < \gamma_2 < \gamma_3,
  \]
  \[
  \varsigma \in (0,1], \;\;\;
  \omega \in \Big(0,\min\big[\half \eta_1, \quarter(1-\eta_2)\big]\Big),\;\;\;
  \kappa_\zeta > \min_{j\in\ii{q}} \epsilon_j^{q+1}
  \tim{and}
  \vartheta_d< \kappa_\zeta.
  \vspace*{-1mm}
  \]
  Choose $\zeta_{d,0}\leq \kappa_\zeta$ and
  set $k=0$ and $i_\zeta=0$.
\item[Step~1: Termination test.]
  Apply the \al{STEP1} algorithm (p.~\pageref{algorithm:STEP1}), resulting in either
  termination, or a model degree $j$ and the associated displacement $d_{k,j}$
  and decrease $\barDT_{f,j}(x_k,d_{k,j})$.
\item[Step~2: Step computation.]
  Apply the \al{STEP2} algorithm (p.~\pageref{algorithm:STEP2}) to compute a
  step $s_k$ such that $\barDT_{f,j}(x_k,s_k)\geq \barDT_{f,j}(x_k,d_{k,j})$.
\item[Step~3: Accept the new iterate.] If
  $\barDT_{f,j}(x_k,s_k)\leq\vartheta_f/\omega$, then terminate with
  $\tilde{x}=x_k$, \status\ = \innoisef, \order\ = $j$, \dlta\ = $\delta_k$
  and \orad = $\max[\delta_k,\|s_k\|]$.\\
  Otherwise, compute $\barf(x_k+s_k)$ ensuring that
  \beqn{Df+-DT}
  |\barf(x_k+s_k)-f(x_k+s_k)| \leq \omega \barDT_{f,j}(x_k,s_k);
  \vspace*{-1mm}
  \eeqn
  and ensure (by setting
  $\barf(x_k)=\barf(x_{k-1}+s_{k-1})$
  or by recomputing $\barf(x_k)$) that 
  \beqn{Df-DT}
  |\barf(x_k)-f(x_k)| \leq \omega \barDT_{f,j}(x_k,s_k).
  \eeqn
  Then compute
  \beqn{rhokdef2}
  \rho_k = \frac{\barf(x_k) - \barf(x_k+s_k)}
                {\barDT_{f,j}(x_k,s_k)}.
  \eeqn
  If $\rho_k \geq \eta_1$, set $x_{k+1} = x_k + s_k$; otherwise set $x_{k+1} = x_k$.
\item[Step~4: Update the trust-region radius.]
  Set
 \[
  \Delta_{k+1} \in \left\{ \begin{array}{ll}
  {}[\gamma_1 \Delta_k, \gamma_2 \Delta_k] & \tim{if} \rho_k < \eta_1,\\
  {}[\gamma_2 \Delta_k, \Delta_k] & \tim{if} \rho_k \in  [\eta_1, \eta_2),\\
  {}[\Delta_k, \min(\Delta_{\max},\gamma_3 \Delta_k)] & \tim{if} \rho_k \geq  \eta_2,
  \end{array}\right.
  \]
  Increment $k$ by one and go to Step~1.
  \end{description}
}

\noindent
We immediately note the condition, at the beginning of Step~3, that
$\barDT_{f,j}(x_k,s_k)>\vartheta_f/\omega$.  This guarantees that the limit
in noise imposed by \req{noise-termination} will not come into play when
computing $\barf(x_k+s_k)$ (and possibly recomputing $\barf(x_k)$).

We also note that, except for
that feature, our specialized \al{STEP1} and \al{STEP2} using the \al{CHECK}
algorithm to handle intrinsic noise on the derivatives, and the relevant
initialization of $\omega$, the
\al{TR$q$EDAN} algorithm is identical to that analyzed in
\cite{CartGoulToin20c}\footnote{\cite{CartGoulToin20c} uses degree-specific
values for $\zeta_d$, but the can be assumed to be identical.}. Again, this
allows us to reuse results in this reference as needed, the first of which
relates the number of iterations of ``successful'' iterations (those where the
new iterate is accepted in Step~3) and ``unsuccessful'' ones.  If, as is
standard, we define 
\[
\calS_k = \{ i \in \iiz{k} \mid x_{i+1}= x_i+s_i\}
= \{ i \in \iiz{k} \mid \rho_i \geq \eta_1 \},
\]
the following useful result is applicable to the \al{TR$q$EDAN} algorithm.

\llem{lemma:SvsU}{\cite[Lemma~3.1]{CartGoulToin20c}
Suppose that the \al{TR$q$EDAN} algorithm is used and that $\Delta_k \geq
\Delta_{\min}$ for some $\Delta_{\min} \in (0,\Delta_0]$. Then, if $k$ is the index of an
iteration before termination,
\beqn{unsucc-neg}
k \leq |\calS_k| \left(1+\frac{\log\gamma_3}{|\log\gamma_2|}\right)
+ \frac{1}{|\log\gamma_2|}
\left|\log\left(\frac{\Delta_{\min}}{\Delta_0}\right)\right|.
\eeqn
}

\subsection{Evaluation complexity and optimality at termination}

Readers with some background in evaluation complexity analysis will not be
surprised by the fact that the complexity of the \al{TR$q$EDAN} algorithm
crucially depends on the decrease that can be achieved on the exact objective
function at successful iterations.  This will in turn depend on the achievable
decrease in inexact values of the objective, which is itself depending on the
decrease $\barDT_{f,j}(x_k,s_k)$ on the inexact model.  Fortunately, we can
again call on the analysis of \cite{CartGoulToin20c} for help, since such
decreases necessarily happen in the \al{TR$q$EDAN} algorithm, before early
termination due to \req{noise-termination} possibly occurs.

\llem{lemma:phihat}{\cite[Lemmas 3.4 and 3.6]{CartGoulToin20c}
Suppose AS.1 and AS.2 hold. At iteration $k$ before termination of the \al{TR$q$EDAN} algorithm
define
\beqn{phihat-def}
\sphi_{f,k} \eqdef \frac{j!\,\,\barDT_{f,j}(x_k,d_{k,j})}{\delta_k^{j}},
\eeqn
where $j$ is the model's degree resulting from \al{STEP1} at iteration $k$. Then,
\beqn{phi-hat-eps}
\sphi_{f,k}\geq \frac{\varsigma\epsmin}{1+\omega},
\eeqn
with $\epsmin=\min_{j\in \{1,\ldots,q\}} \epsilon_j$.
Moreover,
\beqn{DT-phi}
\Delta T_{f,j}(x_k,s_k) \geq \sphi_{f,k}\frac{\delta_k^{j}}{j!}
\tim{ and }
\Delta_k
\geq  \min\left\{\gamma_1\theta,\kappa_r \min_{i\in\iiz{k}}\sphi_{f,i}\right\}
\eeqn
where
\beqn{kappar-def}
\kappa_r
\eqdef \frac{\gamma_1(1-\eta_2)}{4 \max[1,L_f]}
\min\left[\theta,\frac{\Delta_0 \min_{j=1,\ldots,q} \delta_{0,j}^j }
{2q(\max_{j=1,\ldots,q}\|\nabla _x^i f(x_0)\|+\kappa_{\zeta})}\right] \in (0,1).
\eeqn
}

\noindent
Using these results, we may consider the all-important lower bound on the
model decrease at successful iterations.

\llem{lemma:noisy-decrease}{
Suppose that $\vartheta_d >0$ and let $k$ be the index of a successful
iteration of the \al{TR$q$EDAN} algorithm before termination. Then 
\begin{equation}
\label{decr}
\barDT_{f,j}(x_k,s_k)
\geq
\frac{\vartheta_d}{\omega}\varsigma\kappa_{\delta}\epsilon_{\min},
\end{equation}
where
\beqn{kappddef}
\kappa_\delta \eqdef \frac{\kappa_r}{1+\omega}.
\eeqn
}

\proof{
Observe first that, since iteration $k$ is successful, the algorithm must have
reached the end of Step~3 at this iteration, and thus termination did not
occur in Steps 1 or 2. This means in particular, in view of
\req{improve-possible}, that  
\beqn{good-zetas}
\zeta_{d, i_\zeta} > \vartheta_d
\eeqn
for all values of the accuracy threshold $\zeta_{d, i_\zeta}$ encountered during
Steps 1 and 2 of iteration $k$. Moreover Lemma~\ref{lemma:no-term-step1}
applies and \req{barDTd-lower} and \req{good-zetas} imply that
\beqn{4-22}
\barDT_{f,j}(x_k,d_{k,j})
\geq \frac{\zeta_{d,i_\zeta}}{\omega} \sum_{\ell=1}^j\frac{\delta_k^\ell}{\ell!}
\geq \frac{\vartheta_d}{\omega}\delta_k,
\eeqn
again irrespective of the accuracy threshold $\zeta_{d, i_\zeta}$ encountered during
Steps 1 and 2.

We now distinguish two cases, depending on the ratio between $\Delta_k$ and $\theta$.\\
\noindent
{\bf $\bullet$}  Suppose first that $\Delta_k\leq \theta$, (or, equivalently, that
$\delta_k=\Delta_k$).  Then, using \req{step-subproblem} and \req{4-22},  we obtain
that 
\beqn{4-22-b}
\barDT_{f,j}(x_k,s_k)
= \barDT_{f,j}(x_k,d_{k,j})
\geq \frac{\vartheta_d}{\omega}\delta_k
\eeqn
Now, since $\delta_k=\Delta_k \leq \theta$, \req{phi-hat-eps} and the second part
\req{DT-phi} in Lemma~\ref{lemma:phihat} ensure that
$\delta_k \geq \kappa_r \varsigma \epsmin / (1+\omega)$. Substituting this latter bound
in \req{4-22-b} then yields
\[
\barDT_{f,j}(x_k,s_k)
\geq \frac{\vartheta_d\,\kappa_r\,\epsmin}{\omega(1+\omega)}
\]
and \req{decr} follows.\\
\noindent
{\bf $\bullet$}  Suppose now that $\Delta_k> \theta$, (or, equivalently, that
$\delta_k<\Delta_k$).  Then $\delta_k = \theta$. Suppose first that $\|s_k\|\geq
\delta_k=\theta$.  Lemma~\ref{lemma:step2-res} ensures that  \al{STEP2}
terminates with \accuracys being \relative\ and  \req{verif-term-2} holds for
$x=x_k$ and $v=s_k$. As a consequence, using \req{kappar-def}, \req{kappddef} and the fact
that $\varsigma\epsilon_{\min}\leq 1$,
\[
\barDT_{f,j}(x_k,s_k)
\geq \frac{\zeta_{d,i_\zeta}}{\omega} \sum_{\ell=1}^r \frac{\delta_k^\ell}{\ell!}
> \frac{\vartheta_d}{\omega}\sum_{\ell=1}^r \frac{\delta_k^\ell}{\ell!}
\geq \frac{\vartheta_d}{\omega}\theta
\geq \frac{\vartheta_d}{\omega}\varsigma\kappa_\delta \epsilon_{\min},
\]
again implying \req{decr}.  Suppose finally that $\|s_k\|< \delta_k=\theta$.
Then we deduce from  \req{4-22} and \req{step-subproblem} that
\[
\barDT_{f,j}(x_k,s_k)
\geq \barDT_{f,j}(x_k,d_{k,j})
\geq \frac{\vartheta_d}{\omega}\delta_k
= \frac{\vartheta_d}{\omega}\theta
\geq \frac{\vartheta_d}{\omega}\varsigma\kappa_\delta \epsilon_{\min}
\]
and \req{decr} also holds in this last case.
}

\noindent
Of course, this lemma does not allow any useful conlusion if $\vartheta_d=0$,
that is in the noiseless case.  But we can call on the noiseless analysis of
\cite{CartGoulToin20c} to cover this case.

\llem{lemma:noiseless-decrease}{
Suppose that $\vartheta_d=0$. Then, for every $k$ before termination,
\beqn{noiseless-decrease}
\Delta T_j(x_k,s_k) \geq \frac{1}{q!}(\varsigma\kappa_\delta)^{q+1} \epsmin^{q+1},
\eeqn
where $\kappa_\delta$ is defined in \req{kappddef}.
}

\proof{See \cite[Lemma~3.7]{CartGoulToin20c}.  The proof is based on using
  \req{DT-phi} and \req{phi-hat-eps} in Lemma~\ref{lemma:phihat}.
} 

\noindent
We may finally combine Lemmas~\ref{lemma:step1-termination},
\ref{lemma:step2-termination}, \ref{lemma:noisy-decrease} and
\ref{lemma:noiseless-decrease} to derive a upper-bound on the number of
evaluations required by the \al{TR$q$EDAN} algorithm for termination. 

\lthm{theorem:trqedan-complexity}{
Suppose that AS.1--AS.3 hold and define $\epsmin =
\min_{j\in\ii{q}}\epsilon_j$. Then there exists positive constants
$\kappa^A_{{\sf TRqEDAN}}$, $\kappa^B_{{\sf TRqEDAN}}$,
$\kappa^C_{{\sf TRqEDAN}}$, $\kappa^D_{{\sf TRqEDAN}}$,
$\kappa^E_{{\sf TRqEDAN}}$ and $\kappa^S_{{\sf TRqEDAN}}$
such that the \al{TR$q$EDAN} algorithm needs at most
\beqn{number-derivs}
\begin{array}{l}
\kappa^{S}_{{\sf TRqDAN}}\,\bigfrac{f(x_0)-\flow}
         {\max[\vartheta_f,\vartheta_d\epsmin,\epsmin^{q+1}]}
         + \kappa^D_{{\sf TRqDAN}}
        \left|\log\left(\epsmin\right)\right|
+ \kappa^E_{{\sf TRqEDAN}}
\\*[3ex]
\hspace*{6cm}
= \calO\left(\min\left[\vartheta_f^{-1},(\vartheta_d\epsmin)^{-1},\epsmin^{-(q+1)}\right]\right)
\end{array}
\eeqn
evaluations of the (inexact) derivatives $\{\nabla_x^\ell f(x)\}_{\ell=1}^q$, and at most
\beqn{number-f}
\begin{array}{l}
  \kappa^A_{{\sf TRqEDAN}}\bigfrac{f(x_0)-\flow}
        {\max\left[\vartheta_f,\vartheta_d\epsmin,\epsmin^{q+1}\right]}
+ \kappa^B_{{\sf TRqEDAN}}\big|\log(\epsmin)\big|  + \kappa^C_{{\sf TRqEDAN}} \\*[3ex]
\hspace*{6cm}
= \calO\left(\min\left[ \vartheta_f^{-1},(\vartheta_d\epsmin)^{-1},\epsmin^{-(q+1)}\right]\right)
\end{array}
\eeqn
evaluations of $f(x)$ itself to terminate with flags \status, \order, \dlta,
\orad\ and a point $\tilde{x}$ at which
\beqn{optimality0}
\phi_{f,i}^\delta(\tilde{x}) \leq \epsilon_i\frac{\delta^i}{i!}
\tim{for} i\in\ii{j-1}
\eeqn
and
\begin{align}
{\bf \bullet} & \ms
\phi_{f,i}^\delta(x) \leq \epsilon_i \frac{\delta^i}{i!} \tim{for} i \in \iibe{j}{q}\label{optimality1}\\
&\hspace*{7cm} \tim{if \status\ = \asol;}\nonumber\\
{\bf \bullet} & \ms
\phi_{f,j}^\delta(\tilde{x}) \leq \frac{4\vartheta_d}{\gamma_\zeta \omega}\delta\label{optimality2}\\
&\hspace*{7cm} \tim{if \status\ = \innoisephi;}\nonumber\\
{\bf \bullet} & \ms
\phi_{f,j}^\nu(\tilde{x}) \leq \frac{4\vartheta_d}{\gamma_\zeta\omega}
\max\big[\nu,\nu^j\big]\label{optimality3}\\
&\hspace*{7cm} \tim{if \status\ = \innoises,}\nonumber
\end{align}
where  $j$ = \order, $\delta$ = \dlta\ and $\nu$ = \orad.
If, in addition,
\beqn{sk-globmin}
s_k =\arg\max_{\|s_k\|\leq \Delta_k} \barDT_{f,j}(x_k,s)
\eeqn
at iteration $k$ at which termination occurs with \status\ = \innoisef,
then
\beqn{optimality4}
\phi_{f,j}^\nu(\tilde{x}) \leq
\frac{\vartheta_f}{\varsigma}(1+\frac{1}{\omega}).
\eeqn
}

\proof{
We note that the various flag-dependent optimality guarantees
\req{optimality0}--\req{optimality3} are a simple compilation of the results of
Lemmas~\ref{lemma:step1-termination} and \ref{lemma:step2-termination}. To
prove \req{optimality4}, observe that, if termination occurs in Step~3 (as
indicated by \status\ = \innoisef), it
must be because $\barDT_{f,j}(x_k,s_k)\leq \vartheta_f/\omega$.  But
\req{deltak-def} and \req{sk-globmin} imply that
\begin{align*}
\barphi_{f,j}^{\delta_k}(x_k)
&  =  \barDT_{f,j}(x_k,s_k)  \leq \bigfrac{\vartheta_f}{\omega}  \mbox{\hspace*{10mm}if } \|s_k\|\leq \delta_k,\\*[2ex]
\barphi_{f,j}^{\|s_k\|}(x_k)
&\leq \barDT_{f,j}(x_k,s_k)  \leq \bigfrac{\vartheta_f}{\omega}  \mbox{\hspace*{10mm}if } \|s_k\| > \delta_k.
\end{align*}
Moreover, the fact that Step~3 has been reached ensures that termination did
not occur in either Step~1 or Step~2. Thus \req{phi-end-step1} in
Lemma~\ref{lemma:no-term-step1} with the definition \orad\ = $\max[\delta_k,\|s_k\|]$ gives \req{optimality4}.

We now focus on proving \req{number-derivs} and \req{number-f}.  Let $k$ be the
index of a successful iteration before termination. Because \req{Df+-DT} and
\req{Df-DT} both hold at every successful iteration  before termination, we
have that, for each $i\in\calS_k$
\[
f(x_i)-f(x_{i+1})
\geq [\barf(x_i) - \barf(x_{i+1})] - 2 \omega \barDT_{f,j}(x_i,s_i)
\geq (\eta_1 -2 \omega) \barDT_{f,j}(x_i,s_i).
\]
Combining now this inequality with Lemmas~\ref{lemma:noisy-decrease} and
\ref{lemma:noiseless-decrease} we obtain that
\beqn{decr1}
f(x_i)-f(x_{i+1})
\geq (\eta_1 -2 \omega) \max\left[\frac{\vartheta_d}{\omega}\varsigma\kappa_{\delta}\epsilon_{\min},
  \frac{1}{q!}(\varsigma\kappa_\delta)^{q+1} \epsmin^{q+1}\right].
\eeqn
Moreover, the mechanism of Step~3 of the \al{TR$q$EDAN} algorithm implies that
\beqn{decr2}
f(x_i)-f(x_{i+1}) > \frac{\eta_1-2\omega}{\omega}\,\vartheta_f.
\eeqn
From \req{decr1} and \req{decr2}, we thus deduce that
\[
f(x_i)-f(x_{i+1})
\geq (\eta_1 -2 \omega)
  \max\left[\frac{\vartheta_d}{\omega}\varsigma\kappa_{\delta}\epsilon_{\min},
      \frac{1}{q!}(\varsigma\kappa_\delta)^{q+1}\epsmin^{q+1},\frac{\vartheta_f}{\omega}\right]
\eqdef \Delta_f.
\]
Using now the standard ``telescoping sum'' argument and AS.3, we obtain that
\[
f(x_0) - \flow
\geq f(x_0)-f(x_{k+1})
= \sum_{i\in \calS_k} [ f(x_i)- f(x_{i+1})]
\geq |\calS_k| \Delta_f,
\]
so that the total number of successful iterations before termination is 
\beqn{Sk-bound}
|\calS_k|
\leq \frac{f(x_0)-\flow}{\Delta_f}
= \kappa^S_{{\sf TRqEDAN}} \frac{f(x_0)-\flow}
    {\max\left[\vartheta_f,\vartheta_d\epsmin,\epsmin^{q+1}\right]}
\eeqn
where
\[
\kappa^S_{{\sf TRqEDAN}}
\eqdef \frac{1}{(\eta_1-2\omega)}
       \max\left[\frac{1}{\omega},\frac{(\varsigma\kappa_\delta)^{q+1}}{q!}\right]^{-1}.
\]
Now \req{phi-hat-eps}, the second part of \req{DT-phi} and \req{kappddef} imply that
\beqn{rkmin}
\Delta_k\geq \varsigma\kappa_\delta \epsmin,
\eeqn
so that, invoking now Lemma~\ref{lemma:SvsU}, we deduce that the total number of
iterations before termination is bounded above by
\[
n_{\rm it}
\eqdef \frac{f(x_0)-\flow}{\Delta_f}\, \left(1+\frac{\log\gamma_3}{|\log\gamma_2|}\right)
+ \frac{1}{|\log\gamma_2|}
\left|\log\left(\frac{\varsigma\kappa_\delta \epsmin}{\Delta_0}\right)\right|.
\]
Since each iteration of the \al{TR$q$EDAN} algorithm inexactly compute the
objective function's value at most twice (in Step~3),
we obtain that the total number of such evaluations before termination is
bounded above by $2n_{\rm it}$, yielding \req{number-f} with
\[
\kappa^A_{{\sf TRqEDAN}}
\eqdef \frac{2}{\eta_1-2\omega} \min\left[\omega,\frac{q!}{(\varsigma\kappa_\delta)^{q+1}}\right]
       \left(1+\frac{\log\gamma_3}{|\log\gamma_2|}\right),
\]
\[
\kappa^B_{{\sf TRqEDAN}} \eqdef \frac{2}{|\log\gamma_2|}
\tim{ and }
\kappa^C_{{\sf TRqEDAN}} \eqdef \frac{2}{|\log\gamma_2|}
\left|\log\left(\frac{\varsigma\kappa_\delta}{\Delta_0}\right)\right|.
\]

To complete the proof, we need to elaborate on \req{Sk-bound} to derive an
upper bound on the number of derivatives evaluations. While the \al{TR$q$EDAN}
algorithm evaluates $\{\nabla_x^\ell f(x_k)\}_{\ell=1}^j$ at least once in
Step~1, it may need to evaluate the derivatives also when \al{CHECK} returns
\insufficient, and this can happen in the loops between Steps~1.4 and 1.1 in
Algorithm~\ref{algorithm:phij} and between Steps~2.4 and 2.1 in the \al{STEP2}
algorithm. Thus the total number of derivatives' evaluations is given by
$|\calS_k|$ plus the total number of accuracy tightenings (counted by
$i_\zeta$).  The next step is therefore to establish an upper bound on this
latter number. This part of the proof is a variation on that of Theorem~3.8 in
\cite{CartGoulToin20c}, now involving the bounds \req{term-bound-phi} and
\req{step2-bound} but also the additional inequality $\zeta_{d,i_\zeta} \geq
\vartheta_d$ which must hold as long as termination has not
occured.  To summarize the argument, these three bounds ensure a global lower
bound $\zeta_{d,\min}$ on $\zeta_{d,i_\zeta}$, while an upper bound is given
by $\kappa_\zeta$.  Since each tightening proceeds by multiplying the accuracy
threshold by $\gamma_\zeta$, one then deduces that the maximum number of such
tightenings is $\calO\big(|\log(\zeta_{d,\min}/\kappa_\zeta)|\big)$, which then leads to
\req{number-derivs}.  The details are given in appendix.
} 

\noindent
Observe that condition \req{sk-globmin} needs only to be enforced if the bound
\req{optimality4} is desired and when termination occurs with \status\ =
\innoisef. Should \req{optimality4} be of interest, the step
may have to be recomputed in the course of the algorithm to
ensure \req{sk-globmin},  whenever $\barDT_{f,j}(x_k,s_k) <
\vartheta_f/\omega$. Termination is then declared if this inequality still
holds for the new step, or the algorithm is continued otherwise.

The results of Theorem~\ref{theorem:trqedan-complexity} merit some comments.
Firstly, and as expected, we see in the bounds \req{number-derivs} and
\req{number-f} that the total number of evaluations needed for the
\al{TR$q$EDAN} to terminate may be considerably smaller when intrinsic noise
is present ($\vartheta_d >0$ and $\vartheta_f>0$) than in the noiseless
situation ($\vartheta_d=\vartheta_f=0$), in which case we recover the bound in
$\calO(\epsmin^{-(q+1)})+\calO(|\log(\epsmin)|)$ of
\cite{CartGoulToin20c}. More interestingly, we note that, for the intrinsic
noise to be small enough to let the trust-region algoritm run its course
unimpeded, we need that $\vartheta_d = \calO(\epsmin^q)$ and $\vartheta_f =
\calO(\epsmin^{q+1})$.  Since $\vartheta_d$ and $\vartheta_f$ are intrinsic to
the problem, it means that we expect the algorithm to run unimpeded (in the
worst case) only if \beqn{order-observation} \epsmin \gtrsim
\max\left[\vartheta_f^{\frac{1}{q+1}},\vartheta_d^{\frac{1}{q}}\right].  \eeqn
To give an example, suppose that we are applying the \al{TR$q$EDAN} algorithm
to find second-order approximate minimizers on a machine whose machine
precision is $10^{-15}$. This suggest that (in the worst case again), the
algorithm could work as if noise where absent for $\epsmin$ of order $10^{-5}$
and above. Of course, this ignores that some of the deterministic bounds we
have imposed could fail and yet the algorithm could proceed without trouble.
  
We also note that the second term in \req{number-derivs}, which accounts for
the additional evaluations due to inexact but still acceptable evaluations,
now involves a term in $|\log(\vartheta_d/\kappa_\zeta)|$ (the magnitude of the
accuracy range between it initial value and noise) along with the term in
$\log(\epsmin)=\log(\epsmin^q)$ of \cite{CartGoulToin20c}. This is coherent
with our observation \req{order-observation}.

We finally note the difference between the impact of the absolute noise on
the objective function's values ($\vartheta_f$) and that on the derivatives
($\vartheta_d$), the former being significantly more limitative than the
latter. This is reminiscent of similar observations and assumptions
in the stochastic context \cite{BeraCaoSche19,BellGuri21,BlanCartMeniSche19}.

\numsection{Conclusions and Perspectives}
\label{section:conclusion}

We have discussed the evaluation complexity of trust-region algorithms in the
presence of intrinsic noise on function and derivatives values, possibly
causing early termination of the minimization method. We have produced an
evaluation complexity bound which stresses this dependence and relates it to
the complexity bound for the noiseless, albeit inexact, case. In our analysis,
we have priviledged focus and clarity over generality. We have already
mentioned that the noise levels and accuracy thresholds could be made
dependent on the degree of the derivative considered, but other extensions are
indeed possible.  The first is to consider constrained problems, where the
feasible set is convex (or even ``inexpensive'' or ``simple'', see
\cite{BellGuriMoriToin19,CartGoulToin20b,CartGoulToin20a}).  The second is to
replace the Lipschitz continuity required in AS.2 by the weaker H\"{o}lder
continuity (as in
\cite{CartGoulToin17d,CartGoulToin18a,CartGoulToin16,GrapNest17,Nest15}). The
minimization of composite function (using techniques of
\cite{CartGoulToin20b,GratSimoToin20,Nest13}) is another possibility.

Finally, considering ``noise-aware'' stochastic minimization algorithm is also
of interest, and will be reported on in a forthcoming report.
  
{\footnotesize


\begin{thebibliography}{10}

\bibitem{BandScheVice14}
A.~S Bandeira, K.~Scheinberg, and L.~N. Vicente.
\newblock Convergence of trust-region methods based on probabilistic models.
\newblock {\em SIAM Journal on Optimization}, 24(3):1238--1264, 2014.

\bibitem{BellGuri21}
S.~Bellavia and G.~Gurioli.
\newblock Complexity analysis of a stochastic cubic regularisation method under
  inexact gradient evaluations and dynamic hessian accuracy.
\newblock {\em Optimization}, \textmd{(to appear)}, 2021.
\newblock also arXiv:2001.10827.

\bibitem{BellGuriMori21}
S.~Bellavia, G.~Gurioli, and B.~Morini.
\newblock Adaptive cubic regularization methods with dynamic inexact {H}essian
  information and applications to finite-sum minimization.
\newblock {\em IMA Journal of Numerical Analysis}, 41(1):764--799, 2021.

\bibitem{BellGuriMoriToin19}
S.~Bellavia, G.~Gurioli, B.~Morini, and {Ph.}~L. Toint.
\newblock Adaptive regularization algorithms with inexact evaluations for
  nonconvex optimization.
\newblock {\em SIAM Journal on Optimization}, 29(4):2881--2915, 2019.

\bibitem{BellGuriMoriToin20}
S.~Bellavia, G.~Gurioli, B.~Morini, and {Ph.}~L. Toint.
\newblock High-order evaluation complexity of a stochastic adaptive
  regularization algorithm for nonconvex optimization using inexact function
  evaluations and randomly perturbed derivatives.
\newblock arXiv:2005.04639, 2020.

\bibitem{BeraCaoSche19}
A.~Berahas, L.~Cao, and K.~Scheinberg.
\newblock Global convergence rate analysis of a generic line search algorithm
  with noise.
\newblock arXiv:1910.04055, 2019.

\bibitem{BlanCartMeniSche19}
J.~Blanchet, C.~Cartis, M.~Menickelly, and K.~Scheinberg.
\newblock Convergence rate analysis of a stochastic trust region method via
  supermartingales.
\newblock {\em INFORMS Journal on Optimization}, 1(2):92--119, 2019.

\bibitem{Cart91}
R.~G. Carter.
\newblock On the global convergence of trust region methods using inexact
  gradient information.
\newblock {\em SIAM Journal on Numerical Analysis}, 28(1):251--265, 1991.

\bibitem{CartGoulToin17d}
C.~Cartis, N.~I.~M. Gould, and Ph.~L. Toint.
\newblock Worst-case evaluation complexity of regularization methods for smooth
  unconstrained optimization using {H}\"older continuous gradients.
\newblock {\em Optimization Methods and Software}, 6(6):1273--1298, 2017.

\bibitem{CartGoulToin18a}
C.~Cartis, N.~I.~M. Gould, and Ph.~L. Toint.
\newblock Worst-case evaluation complexity and optimality of second-order
  methods for nonconvex smooth optimization.
\newblock In B.~Sirakov, P.~{de Souza}, and M.~Viana, editors, {\em Invited
  Lectures, Proceedings of the 2018 International Conference of Mathematicians
  (ICM 2018), vol. 4, Rio de Janeiro}, pages 3729--3768. World Scientific
  Publishing Co Pte Ltd, 2018.

\bibitem{CartGoulToin16}
C.~Cartis, N.~I.~M. Gould, and Ph.~L. Toint.
\newblock Universal regularization methods -- varying the power, the smoothness
  and the accuracy.
\newblock {\em SIAM Journal on Optimization}, 29(1):595–--615, 2019.

\bibitem{CartGoulToin20b}
C.~Cartis, N.~I.~M. Gould, and Ph.~L. Toint.
\newblock Sharp worst-case evaluation complexity bounds for arbitrary-order
  nonconvex optimization with inexpensive constraints.
\newblock {\em SIAM Journal on Optimization}, 30(1):513--541, 2020.

\bibitem{CartGoulToin20a}
C.~Cartis, N.~I.~M. Gould, and Ph.~L. Toint.
\newblock Strong evaluation complexity bounds for arbitrary-order optimization
  of nonconvex nonsmooth composite functions.
\newblock arXiv:2001.10802, 2020.

\bibitem{CartGoulToin20c}
C.~Cartis, N.~I.~M. Gould, and Ph.~L. Toint.
\newblock Strong evaluation complexity of an inexact trust-region algorithm for
  arbitrary-order unconstrained nonconvex optimization.
\newblock arXiv:2011.00854, 2020.

\bibitem{ChenMeniSche18}
R.~Chen, M.~Menickelly, and K.~Scheinberg.
\newblock Stochastic optimization using a trust-region method and random
  models.
\newblock {\em Mathematical Programming, Series~A}, 169(2):447--487, 2018.

\bibitem{ConnGoulToin00}
A.~R. Conn, N.~I.~M. Gould, and Ph.~L. Toint.
\newblock {\em Trust-Region Methods}.
\newblock MPS-SIAM Series on Optimization. SIAM, Philadelphia, USA, 2000.

\bibitem{deKlLaur19}
E.~de~Klerk and M.~Laurent.
\newblock Worst-case examples for {L}asserre's measure-based hierarchy for
  polynomial optimization on the hypercube.
\newblock {\em Mathematics of Operations Research}, 45(1):86--98, 2019.

\bibitem{deKlLaur20}
E.~de~Klerk and M.~Laurent.
\newblock Convergence analysis of a {L}asserre hierarchy of upper bounds for
  polynomial minimization on the sphere.
\newblock {\em Mathematical Programming}, \textmd{(to appear)}, 2020.

\bibitem{GrapNest17}
G.~N. Grapiglia and {Yu}. Nesterov.
\newblock Regularized {N}ewton methods for minimizing functions with {H}\"older
  continuous {H}essians.
\newblock {\em SIAM Journal on Optimization}, 27:478--506, 2017.

\bibitem{GratSartToin08}
S.~Gratton, A.~Sartenaer, and Ph.~L. Toint.
\newblock Recursive trust-region methods for multiscale nonlinear optimization.
\newblock {\em SIAM Journal on Optimization}, 19(1):414--444, 2008.

\bibitem{GratSimoToin20}
S.~Gratton, E.~Simon, and Ph.~L. Toint.
\newblock An algorithm for the minimization of nonsmooth nonconvex functions
  using inexact evaluations and its worst-case complexity.
\newblock {\em Mathematical Programming, Series~A}, \textmd{(to appear)}, 2021.

\bibitem{GratToin20}
S.~Gratton and Ph.~L. Toint.
\newblock A note on solving nonlinear optimization problems in variable
  precision.
\newblock {\em Computational Optimization and Applications}, 76(3):917--933,
  2020.

\bibitem{High17}
N.~J. Higham.
\newblock The rise of multiprecision computations.
\newblock Talk at SAMSI 2017, April 2017.
\newblock {\tt https://bit.ly/higham-samsi17}.

\bibitem{Nest13}
{Yu}. Nesterov.
\newblock Gradient mehods for minimizing composite objective functions.
\newblock {\em Mathematical Programming, Series~A}, 140(1):125--161, 2013.

\bibitem{Nest15}
{Yu}. Nesterov.
\newblock Universal gradient methods for convex optimization problems.
\newblock {\em Mathematical Programming, Series~A}, 152(1–2):381--–404,
  2015.

\bibitem{PaquSche20}
C.~Paquette and K.~Scheinberg.
\newblock A stochastic line search method with convergence rate analysis.
\newblock {\em SIAM Journal on Optimization}, 30(1):349--376, 2020.

\bibitem{SlotLaur20}
L.~Slot and M.~Laurent.
\newblock Improved convergence analysis of {L}asserre’s measure-based upper
  bounds for polynomial minimization on compact sets.
\newblock {\em Mathematical Programming}, \textmd{(to appear)}, 2020.

\bibitem{XuRoosMaho20}
P.~Xu, F.~Roosta-Khorasani, and M.~W. Mahoney.
\newblock {N}ewton-type methods for non-convex optimization under inexact
  {H}essian information.
\newblock {\em Mathematical Programming, Series~A}, 184((1-2)):35--70, 2020.

\bibitem{YaoXuRoosMaho20}
Z.~Yao, P.~Xu, F.~Roosta-Khorasani, and M.~W. Mahoney.
\newblock Inexact non-convex {N}ewton-type methods, 2020.

\bibitem{Yuan15}
Y.~Yuan.
\newblock Recent advances in trust region algorithms.
\newblock {\em Mathematical Programming, Series~A}, 151(1):249--281, 2015.

\end{thebibliography}

}

\appendix

\appnumsection{Details of the proof of Theorem~\ref{theorem:trqedan-complexity}}

We follow the argument of \cite[proof of Theorem~3.8]{CartGoulToin20c},
(adapting the bounds to the new context), and derive an upper bound on the
number of derivatives' evaluations. This requires counting the number of
additional derivative evaluations caused by successive tightening of the
accuracy threshold $\zeta_{d,i_\zeta}$. Observe that repeated evaluations at a
given iterate $x_k$ are only needed when the current value of this threshold
is smaller than used previously at the same iterate $x_k$. The
$\{\zeta_{d,i_\zeta}\}$ are, by construction, linearly decreasing with rate
$\gamma_\zeta$,  Indeed, $\zeta_{d,i_\zeta}$ is initialised to
$\zeta_{d,0}\leq \kappa_\zeta$
in Step~0 of the \al{TR$q$DAN} algorithm,  decreased each time by a factor
$\gamma_\zeta$ in \req{improve-accuracy} in the \al{CHECK} invoked in Step~1.2
of Algorithm~\ref{algorithm:phij}, down to the value $\zeta_{d,i_\zeta}$ which
is then passed to Step~2, and possibly decreased there further in
\req{improve-accuracy} in the \al{CHECK} invoked in Step~2.2 of the \al{STEP2}
algorithm again by successive multiplication by $\gamma_\zeta$. We now use
\req{term-bound-phi} in Lemma~\ref{lemma:phi-res} and \req{step2-bound} in
Lemma~\ref{lemma:step2-res} to deduce that, even in the absence of noise,
$\zeta_{d,i_\zeta}$ will not be reduced below the value
\beqn{minval}
\min\left[ \frac{\omega}{4}\,\varsigma\,\epsilon_j\,\frac{\delta_k^{j-1}}{j!},
       \frac{\omega}{8(1+\omega)\max[1,\Delta_{\max}^j]}\,\epsilon_j\,\frac{\delta_k^j}{j!}\right]
\geq \frac{\varsigma\,\omega}{8(1+\omega)\max[1,\Delta_{\max}^j]} \,\epsilon_j\,\frac{\delta_k^j}{j!}
\eeqn
at iteration $k$. Now define
\[
\kappa_{\rm acc}
\eqdef \frac{\varsigma\omega(\varsigma\kappa_{\delta})^q}{8(1+\omega)\max[1,\Delta_{\max}^j]}
\leq  \frac{\varsigma\omega}{8(1+\omega)\max[1,\Delta_{\max}^j]} \,\frac{(\varsigma\kappa_{\delta})^j}{j!}
\]
so that \req{rkmin} implies that
\[
\kappa_{\rm acc}\epsmin^{q+1}
\leq \frac{\varsigma\omega\,\epsilon_j}{8(1+\omega)\max[1,\Delta_{\max}^j]}
\,\frac{\delta_k^j}{j!}.
\]
We also note that conditions \req{improve-possible} and
\req{improve-accuracy} in the \al{CHECK} algorithm impose that any reduced
value of $\zeta_{d,i_\zeta}$ (before termination) must satisfy the bound
$\zeta_{d,i_\zeta} \geq \vartheta_d$. Hence the bound \req{minval} can be
strengthened to be
\[
\max\left[\vartheta_d, \kappa_{\rm acc}\epsmin^{q+1}\right].
\]
Thus no further reduction of the 
$\zeta_{d,i_\zeta}$, and hence no further approximation of
$\{\overline{\nabla_x^j f}(x_k)\}_{j=1}^q$,  can possibly occur in any
iteration once the largest initial absolute error $\zeta_{d,0}$ has been
reduced by successive multiplications by $\gamma_\zeta$ sufficiently to ensure
that
\beqn{zetamax}
\gamma_\zeta^{i_\zeta}\zeta_{d,0}
\leq \gamma_\zeta^{i_\zeta}\kappa_\zeta
\leq \max[\vartheta_d,\kappa_{\rm acc} \epsmin^{q+1}],
\eeqn
the second inequality being equivalent to asking
\beqn{izeta-bound}
i_\zeta \log(\gamma_\zeta) 
\leq \max\left[\log(\vartheta_d), (q+1)\log\left(\epsmin\right) +
  \log(\kappa_{\rm acc})\right]
- \log\left(\kappa_\zeta\right),
\eeqn
where the right-hand side is negative because of the inequalities
$\kappa_{\rm acc}<1$ and $\max[\epsmin^{q+1},\vartheta_d] \leq \kappa_\zeta$ (imposed in the
initialization step of the \al{TRqEDAN} algorithm). 
We now recall that Step~1 of this  algorithm is only used (and
derivatives evaluated) after successful iterations. As a consequence, we
deduce that the number of evaluations of the derivatives of the objective
function that occur during the course of the \al{TR$p$DAN} algorithm before
termination is at most
\beqn{evald-trqda}
|\calS_k| + i_{\zeta,\max},
\eeqn
i.e., the number iterations in \req{Sk-bound} plus
\[
\begin{array}{lcl}
i_{\zeta,\max}
& \!\! \eqdef \!\! &
\left\lfloor 
\bigfrac{1}{\log(\gamma_\zeta)}
\max\left\{
  \log\left(\bigfrac{\vartheta_d}{\zeta_{d,0}}\right),
  (q+1)\log\left(\epsmin \right)+\log\left(\bigfrac{\kappa_{\rm acc}}{\zeta_{d,0}}\right)
  \right\}\right\rfloor \\*[2ex]
& \!\!< \!\! &
  \bigfrac{1}{|\log(\gamma_\zeta)|}
  \left\{
    \left|\log\left(\bigfrac{\vartheta_d}{\zeta_{d,0}}\right)\right| +
    \,(q+1)\left|\log\left(\epsmin\right)\right|
  + \left|\log\left(\bigfrac{\kappa_{\rm acc}}{\zeta_{d,0}}\right)\right|\right\}+1,
\end{array}
\]
the largest value of $i_\zeta$ that ensures \req{izeta-bound}.
Adding one for the final evaluation at termination, this leads to the desired
evaluation bound \req{number-derivs} with the coefficients
\[
\kappa^D_{{\sf TRqEDAN}} \eqdef \frac{q+1}{|\log\gamma_\zeta|}
\tim{and}
\kappa^E_{{\sf TRqEDAN}} \eqdef
\bigfrac{1}{|\log(\gamma_\zeta)|}
\left\{
    \left|\log\left(\bigfrac{\kappa_{\rm acc}}{\zeta_{d,0}}\right)\right|
    + \left|\log\left(\bigfrac{\vartheta_d}{\zeta_{d,0}}\right)\right|
    \right\}
    +2.
\]

\end{document}